\documentclass[a4paper,reqno,oneside,11pt]{amsart}
\usepackage{fullpage}
\usepackage{setspace}
\usepackage{amsmath, amssymb, amsfonts, amsthm, mathtools, mathrsfs}
\usepackage{tikz}
\usepackage{subcaption}
\usepackage{verbatim,enumitem,graphics}
\usepackage{xcolor}
\usepackage[backref=page]{hyperref}
\hypersetup{%
    colorlinks,
    linkcolor={red!50!black},
    citecolor={green!50!black},
    urlcolor={blue!50!black}
}
\usepackage{lineno}

\allowdisplaybreaks

\theoremstyle{plain}
\newtheorem{theorem}{Theorem}[section]
\newtheorem*{theorem*}{Theorem}
\newtheorem{lemma}[theorem]{Lemma}
\newtheorem{proposition}[theorem]{Proposition}

\newtheorem{corollary}[theorem]{Corollary}

\theoremstyle{definition}
\newtheorem{claim}[theorem]{Claim}
\newtheorem{claim*}{Claim}

\newtheorem*{remark*}{Remark}
\newtheorem{definition}[theorem]{Definition}

\newcommand{\scC}{\mathscr{C}}

\newcommand{\scF}{\mathscr{F}}

\newcommand{\scH}{\mathscr{H}}

\newcommand{\scP}{\mathscr{P}}

\newcommand{\scX}{\mathscr{X}}

\newcommand{\bbE}{\mathbb{E}}

\newcommand{\bbP}{\mathbb{P}}

\newcommand{\N}{\mathbb{N}}

\newcommand{\F}{\mathbb{F}}

\newcommand{\fF}{\mathfrak{F}}

\def\d{\delta}

\def\l{\ell}
\def\e{\epsilon}

\begin{document}

\setstretch{1.27}

\title{On some off-diagonal hypergraph Ramsey problems}

\author{Eion Mulrenin}
\address{Department of Mathematics, Emory University, Atlanta, Georgia, USA}
\email{eion.mulrenin@emory.edu}

\author{Cosmin Pohoata}
\address{Department of Mathematics, Emory University, Atlanta, Georgia, USA}
\email{cosmin.pohoata@emory.edu}

\begin{abstract}
    We give a short proof of a theorem of Fox and He that the $3$-uniform Ramsey number $r_3(K_4^-, K_n)$ is at least $2^{\Omega(n \log n)}$.
    The construction exploits a new connection with a directed graph recently used by Brada\v{c} in his work on off-diagonal graph Ramsey numbers.
    Moreover, it generalizes naturally to give sharp bounds for several extensions of this problem in higher uniformity, answering questions raised by Mubayi and Suk.
\end{abstract}

\maketitle


\section{Introduction}
\label{sec: intro}

Ramsey's theorem~\cite{Ramsey} is the statement that for two $k$-uniform hypergraphs $G$ and $H$ and $N$ sufficiently large, any red-blue coloring of the edges of the complete $k$-uniform hypergraph $K_N^{(k)}$ will contain either a red copy of $G$ or a blue copy of $H$.
The quantitative aspects of the least $N$ with this property, which is denoted by $r_k(G, H)$, lay the foundation for many major open problems in extremal combinatorics and beyond.

One regime which has been a central object of study for nearly a century is the {\it off-diagonal Ramsey numbers}, where one thinks of $G$ as some fixed $k$-graph and $H = K_n^{(k)}$ with $n$ going to infinity.
The graph case $k=2$ has recently seen major progress for fixed cliques $G=K_s$: following a breakthrough by Brada\v{c}~\cite{Brad26+}, who nearly obtained the correct power of $n$, it was proved by an internal model at OpenAI~\cite{OpenAI5,OpenAI} that
\[
    r_2(s,n) = \frac{n^{s-1}}{(\log n)^{s-2+o(1)}}
\]
for every fixed $s \geq 5$ as $n\to\infty$, asymptotically matching the logarithmic exponent in the classical upper bound of Ajtai, Koml\'os, and Szemer\'edi~\cite{AKS80}.
Lin and Zhang~\cite{LZ26+} subsequently obtained the same conclusion for $s=4$, while Kim~\cite{Kim} proved the stronger estimate $r_2(3,n)=\Theta(n^2/\log n)$ in 1995.
The central object behind the constructions of Brada\v{c} and OpenAI is a directed graph $D_q$ arising from incidences in finite projective spaces, which had been used before in Ramsey theory (see, e.g.,~\cite{KPR10}) and is well-known in the finite geometry literature~\cite{BMM20, Ott, Smith}.
In this paper, we will show how this same digraph may be used to give constructions for several hypergraph Ramsey problems, that is, when the uniformity $k$ is at least~$3$. Throughout the paper, $k$ is fixed and implicit constants may depend on $k$. When choosing integer parameters, we suppress floor signs and round down where necessary.

The first of these concerns the off-diagonal Ramsey number $r_3(K_4^-, K_n)$, where $K_4^-$ denotes the hypergraph obtained by removing an edge from a copy of $K_4$, or equivalently, three edges on four vertices\footnote{When the uniformity is understood, e.g., here where the subscript $3$ on $r_3(K_4^-, K_n)$ indicates that the uniformity is $3$, we will drop the superscript on cliques and write $K_n$ for $K_n^{(3)}$.}.
If we instead ask for only two red edges on four vertices, Erd\H{o}s and Hajnal~\cite{EH72}observed that the corresponding Ramsey number is polynomial in $n$; in fact, the exact order of magnitude $\Theta(n^2/\log n)$ is known~\cite{KMV14, RS94}, with similarly sharp results holding in uniformity $k$.
Perhaps surprisingly, then, Erd\H{o}s and Hajnal proved in the same paper~\cite{EH72} the exponential-type bounds
\[
    2^{\Omega(n)} \leq r_3(K_4^-, K_n) \leq 2^{O(n \log n)}
\]
for three red edges on four vertices,
and conjectured that the upper bound was sharp up to the constant in the exponent.
This was resolved several years ago by Fox and He~\cite{FH21}.

\begin{theorem}[Fox--He, 2019~\cite{FH21}]
\label{thm: k4-}
    $r_3(K_4^-, K_n) \geq 2^{\Omega(n \log n)}$.
\end{theorem}

The lower bound of Erd\H{o}s and Hajnal follows from considering a random tournament on $N$ vertices and coloring red the triples forming a cyclic triangle, i.e., triples $xyz$ in which $x \to y \to z \to x$ in the orientation of edges, and all other triples blue.
One can check that there are no red copies of $K_4^-$, and a calculation of the expected number of sets of size $\Omega(\log N)$ with all blue edges gives the result.
Fox and He then boosted this idea with a rather complicated entropy argument.
Our first result is a short proof of this result, which also builds on the idea behind the Erd\H{o}s--Hajnal construction by exploiting the aforementioned directed graph $D_q$ in which all triangles are cyclic.

Our remaining results concern generalizations of this problem, each of which goes back to the following, introduced by Erd\H{o}s and Hajnal~\cite{EH72} in 1972.

\begin{definition}
    For integers $2 \leq k \leq s \leq n$ and $2 \leq t \leq \binom{s}{k}$, let $r_k(s, t; n)$ denote the least integer $N$ such that every red-blue coloring of $E(K_N^{(k)})$ contains either a set of $s$ vertices spanning at least $t$ red edges or a blue copy of $K_n^{(k)}$.
\end{definition}
Note that $r_k(s, \binom{s}{k}; n) = r_k(s,n)$ recovers the usual Ramsey numbers, while $r_3(4, 3; n) = r_3(K_4^-, K_n)$.
Erd\H{o}s and Hajnal generalized their lower bound to higher uniformities, showing that
\begin{equation}
\label{eq: erdos-hajnal-original-bound}
    r_k(k+1, 2; n) \leq O(n^{k-1})
    \qquad
    \text{and}
    \qquad
    r_k(k+1, 3; n) \geq 2^{\Omega(n)},
\end{equation}
and conjectured more generally that as $t$ increases from $1$ to $\binom{s}{k}$, there should be a precise choice $t_1 = h_1^{(k)}(s)$ at which $r_k(s,t; n)$ jumps from polynomial to exponential in a power of $n$, a second choice $t_2 = h_2^{(k)}(s)$ at which it jumps from exponential to doubly exponential in a power of $n$, and so on, up to a choice $t_{k-2} = h_{k-2}^{(k)}(s)$ at which it jumps from a $(k-2)$- to a $(k-1)$-fold iterated exponential in a power of $n$.
The first threshold $h_1^{(k)}(s)$ is now known for all $s>k\geq3$. Mubayi and Razborov~\cite{MR21} proved the conjectured polynomial-to-exponential transition for every uniformity $k\geq4$ and every $s>k$. In uniformity $3$, Conlon, Fox, and Sudakov~\cite{CFS12} established it for infinitely many values of $s$, and Ascoli, He, and Yu~\cite{AHY25+} settled the remaining cases.

Our next result concerns $r_k(k+1, 3; n)$, which recovers $r_3(K_4^-, K_n)$ when $k=3$.
Mubayi and Suk~\cite{MS20a} improved the lower bound in~\eqref{eq: erdos-hajnal-original-bound} and proved a similar upper bound, showing that
\[
    2^{\Omega(n^{k-2})} \leq r_k(k+1, 3; n) \leq 2^{O(n^{k-2} \log n)}.
\]
Our proof of Theorem~\ref{thm: k4-} directly generalizes to show that the upper bound is again tight up to the constant in the exponent.

\begin{theorem}
\label{thm: 3-edges-k-unif}
    For every fixed $k\geq3$, we have $r_k(k+1, 3; n) \geq 2^{\Omega_k(n^{k-2} \log n)}$.
\end{theorem}

Another generalization of $r_3(K_4^-, K_n)$, which goes back to the 1972 paper of Erd\H{o}s and Hajnal~\cite{EH72}, is the Ramsey number of the {\it half-graph} $B^{(k)}$ versus $K_n^{(k)}$.
The $k$-half-graph $B^{(k)}$ is a $k$-uniform hypergraph with $2k-2$ vertices and {$k$ edges} defined as follows: its vertex set is a disjoint union $S \cup T$, where $|S|=|T|=k-1$, and its edge set is {$\{S \cup \{t\}: t \in T\} \cup \{T \cup \{s_*\}\}$} for some distinguished element $s_* \in S$, i.e., all $k$-sets containing $S$ and exactly one $k$-set containing $T$.
Observe that $K_4^-$ is just the 3-half-graph $B^{(3)}$.
Generalizing the bounds $2^{\Omega(n)} \leq r_3(K_4^-, K_n) \leq (n!)^2$ of Erd\H{o}s and Hajnal~\cite{EH72}, Mubayi and Suk~\cite{MS18} showed that
\[
    2^{\Omega(n)} \leq r_k(B^{(k)}, K_n^{(k)}) \leq (n!)^{k-1},
\]
and posed the problem of improving the lower bound~\cite[Problem 8.1]{MS20b}.
Our third result shows that once again, the upper bound is tight up to the constant in the exponent.

\begin{theorem}
\label{thm: half-graph}
    For every fixed $k\geq3$, we have $r_k(B^{(k)}, K_n^{(k)}) \geq 2^{\Omega_k(n \log n)}$.
\end{theorem}

Finally, we observe that our new bound in Theorem~\ref{thm: 3-edges-k-unif} may be combined with results of Mubayi and Suk~\cite{MS20a} to obtain better bounds on $r_k(k+1, t; n)$ for values of $t \geq 4$.
To state these bounds, recall that the {\it tower function} $T_j(x)$ is defined recursively by $T_1(x) = x$ and $T_{j+1}(x) = 2^{T_j(x)}$ for $j \geq 1$.
Mubayi and Suk, using the stepping-up method of Erd\H{o}s and Hajnal~\cite{EHR65}, showed that for $k \geq 6$ and $4 \leq t \leq k-2$,
\[
    T_{t-1}(O(n^{k-t+1} \log n)) \geq 
    r_k(k+1, t; n) \geq 
    \begin{cases}
        T_{t-1}(c n^{k-t+1}), \quad k-t \text{ even};\\
        T_{t-1}(c n^{\frac{k-t+1}{2}}), \quad k-t \text{ odd}.
    \end{cases}
\]

Our final result adds logarithms to the tops of the towers in this lower bound, matching the upper bound up to the constant when $k-t$ is even.

\begin{theorem}
\label{thm: t-edges-k-unif}
    For all $k \geq 6$ and $4 \leq t \leq k-2$, there exists a constant $c = c(k, t)>0$ such that the following holds for all sufficiently large $n$:
    \begin{equation}
        r_k(k+1, t; n) \geq 
        \begin{cases}
            T_{t-1}(c n^{k-t+1} \log n), \quad k-t \text{ even};\\
            T_{t-1}(c n^{\frac{k-t+1}{2}} \log n), \quad k-t \text{ odd}.
        \end{cases}
    \end{equation}  
    In particular, when $k-t$ is even, we have
    \[
        r_k(k+1,t;n) = T_{t-1}(\Theta(n^{k-t+1} \log n)).
    \]
\end{theorem}

\subsection*{Organization}
The remainder of the paper is organized as follows: in Section~\ref{sec: K4-}, we give a proof of the bound $r_3(K_4^-, K_n) \geq 2^{\Omega(n \log n)}$; in Section~\ref{sec: ind-number}, we prove a general lemma which we will use to bound the size of the largest blue clique in the proofs of Theorems~\ref{thm: 3-edges-k-unif} and~\ref{thm: half-graph}; in Section~\ref{sec: proofs-k-unif}, we prove these theorems; and in Section~\ref{sec: proofs-larger-t}, we prove Theorem~\ref{thm: t-edges-k-unif}.


\section{The Ramsey number \texorpdfstring{$r_3(K_4^-, K_n)$}{r3(K4-, Kn)}}
\label{sec: K4-}

In this section, we present a proof of Theorem~\ref{thm: k4-}, i.e., that $r_3(K_4^-, K_n) \geq 2^{\Omega(n \log n)}$.

To begin, let $q$ be a prime power, let $\mathrm{PG}(2,q)$ denote the projective plane of order $q$ with points $P_q$ and lines $L_q$, and consider the set
\[
    \fF_q = \{(p, \l): p \in \l\} \subseteq P_q \times L_q
\]
of flags in $\mathrm{PG}(2,q)$.
Since there are $q^2+q+1$ points and each point lies on $q+1$ lines, we have
\[
    M := |\fF_q| = (q^2+q+1) \cdot (q+1).
\]
We will use the fact that 
\begin{equation}
\label{eq: size-M}
    q^3 \leq M \leq 3q^3
\end{equation}
for all prime powers $q$ in some of the calculations below.

Our starting point for the construction will be the directed graph $D_q$, recently introduced by Brada\v{c}, which has vertex set $\fF_q$ and an arc $(p,\l) \to (p',\l')$ if and only if
\[
    p \in \l' \qquad \text{ but } \qquad p' \notin \l.
\]
Geometrically, there is an arc from $(p,\l)$ to $(p', \l')$ exactly when $p \neq p'$, $\l \neq \l'$ and $\l'$ is the unique line determined by $p$ and $p'$.
Note that every vertex of $D_q$ has exactly $q^2$ in-neighbors and $q^2$ out-neighbors, since every point has $q+1$ lines passing through it and every line passes through $q+1$ points.

The following proposition, which we will use in all of the constructions in this paper, is well-known. For a proof, see, e.g.,~\cite[Lemma 2.11]{Brad26+}; for more details and discussion, see~\cite{BMM20}.

\begin{proposition}
\label{prp: Dq-no-TT}
    For all prime powers $q$, the directed graph $D_q$ contains no transitive triangle.
    That is, there do not exist $x_1, x_2, x_3 \in V(D_q)$ with $x_i \to x_j$ for all $i < j \in [3]$.
\end{proposition}

Since $D_q$ contains no transitive triangles, all triangles in the digraph are cyclic.
It turns out that these cyclic triangles are rather robust in $D_q$, in the following sense.

\begin{proposition}
\label{prp: cyclic-triangles}
    If subsets $X,Y,Z\subseteq\fF_q$ satisfy
    \[
        |X||Y|,\ |Y||Z|,\ |X||Z|\geq 5000q^5,
    \]
    then there exist $x\in X$, $y\in Y$, and $z\in Z$ such that
    $x\to y\to z\to x$ in $D_q$.
\end{proposition}

\begin{proof}
    Given a flag $w\in\fF_q$ and a set $W\subseteq\fF_q$, write
    $d^+(w,W)$ and $d^-(w,W)$ for the out- and in-degrees of $w$
    into $W$ in $D_q$, respectively. We first claim that, for a uniformly random
    flag $w\in\fF_q$,
    \begin{equation}
    \label{eq: bradac-variance}
        \mathrm{Var}(d^+(w,W)),\ \mathrm{Var}(d^-(w,W))\leq 10|W|.
    \end{equation}
    Indeed, since $D_q$ is $q^2$-regular,
    \[
        \bbE[d^+(w,W)]=\frac{q^2|W|}{M}.
    \]
    By counting common in-neighbors of ordered pairs of flags in $W$, we get
    \[
        \sum_{w\in\fF_q}d^+(w,W)^2
        \leq q^2|W|+q^3|W|+(q-1)|W|^2
        \leq q|W|^2+2q^3|W|:
    \]
    here we used that two distinct flags on the same line have
    $q(q-1)$ common in-neighbors, while two flags on different lines
    have at most $q-1$ common in-neighbors; and moreover, each flag
    shares its line with at most $q$ other flags. Consequently,
    \begin{align*}
        \mathrm{Var}(d^+(w,W))
        &\leq \frac{q|W|^2+2q^3|W|}{M}-\frac{q^4|W|^2}{M^2}\\
        &=\frac{(qM-q^4)|W|^2}{M^2}+\frac{2q^3|W|}{M}
        \leq 10|W|,
    \end{align*}
    where the last inequality uses $|W|\leq M$ and
    $qM-q^4=2q^3+2q^2+q$. The same bound for $d^-(w,W)$
    follows by taking a polarity of $\mathrm{PG}(2,q)$, which is an automorphism of the underlying graph of $D_q$ that reverses the directions of the arcs,
    proving~\eqref{eq: bradac-variance}.

    Now, since $\bbE[d^+(w,X)]=q^2|X|/M$, Chebyshev's inequality
    and~\eqref{eq: bradac-variance} give
    \[
        \#\left\{w\in\fF_q:
        d^+(w,X)<\frac{q^2|X|}{2M}\right\}
        \leq\frac{40M^3}{q^4|X|}
        \leq\frac{1080q^5}{|X|}<\frac{|Z|}{4}.
    \]
    Similarly,
    \[
        \#\left\{w\in\fF_q:
        d^-(w,Y)<\frac{q^2|Y|}{2M}\right\}<\frac{|Z|}{4}.
    \]
    Thus there exists $z\in Z$ with
    \[
        d^+(z,X)\geq\frac{q^2|X|}{2M},\qquad
        d^-(z,Y)\geq\frac{q^2|Y|}{2M}.
    \]
    Set $A=N^+(z)\cap X$ and $B=N^-(z)\cap Y$. Then
    \[
        |A||B|\geq\frac{q^4|X||Y|}{4M^2}\geq4q^3.
    \]
    We claim that there is an arc from $A$ to $B$. To see this,
    write $z=(p,\l)$ and regard $\l$ as the line at infinity
    in $\mathrm{PG}(2,q)$. The out-neighbors of $z$ correspond
    bijectively to the affine points of $\mathrm{AG}(2,q)$, while
    its in-neighbors correspond to the affine lines with directions
    different from $p$. Under these correspondences, an arc
    from an out-neighbor to an in-neighbor is exactly a point-line
    incidence in $\mathrm{AG}(2,q)$. Hence Vinh's incidence bound~\cite[Theorem 3]{Vinh} yields\footnote{Vinh only states the upper bound $e(A,B)\leq\frac{|A||B|}{q}+\sqrt{q|A||B|}$, but since his proof uses the expander-mixing lemma, the lower bound also follows.}
    \[
        e(A,B)\geq\frac{|A||B|}{q}-\sqrt{q|A||B|}
        \geq\frac{|A||B|}{2q}>0,
    \]
    where $e(A,B)$ denotes the number of arcs from $A$ to $B$.
    Choosing an arc $x\to y$ with $x\in A$ and $y\in B$,
    we obtain $x\to y\to z\to x$, as desired.
\end{proof}

\begin{remark*}
    The $q^5$ threshold in Proposition~\ref{prp: cyclic-triangles} is sharp up to an absolute constant, even when $q$ is prime.
    To see this, let $p$ be a sufficiently large prime and put $r=\sqrt p/4$.
    Using integer representatives for elements of $\F_p$, take
    \[
        P=\{(x,y):1\leq x\leq r,\ 0\leq y\leq p/4\}
    \]
    and the family of affine lines
    \[
        \mathcal L=\{y=mx+b:1\leq m\leq r,\
        p/2\leq b\leq3p/4\}.
    \]
    Both $P$ and $\mathcal L$ have size $\Theta(p^{3/2})$, but there are no incidences between them: for these integer representatives,
    \[
        \frac p2\leq mx+b\leq\frac{3p}{4}+r^2
        \leq\frac{13p}{16}<p,
    \]
    whereas every point of $P$ has second coordinate at most $p/4$.
    Identifying the lines of $\mathcal L$ with their projective completions, set
    \[
        X=Z=\{(u,\ell)\in\fF_p:u\in P\},\qquad
        Y=\{(u,\ell)\in\fF_p:\ell\in\mathcal L\}.
    \]
    Since each point and each line belongs to $p+1$ flags, all three sets have size $\Theta(p^{5/2})$.
    However, there are no arcs from $X$ to $Y$, and hence no cyclic triangles $x\to y\to z\to x$ with $x\in X$, $y\in Y$, and $z\in Z$.
    Thus the pairwise-product threshold cannot be replaced by $o(p^5)$, even over prime fields.
\end{remark*}

Proposition~\ref{prp: cyclic-triangles} yields the following corollaries which we will need below.

\begin{corollary}
\label{cor: cyclic-triangles}
    Let $q$ be a sufficiently large prime power, and let
    $X_1,X_2,X_3\subseteq\fF_q$.
    If $|X_i|\geq M^{6/7}$ for each $i$, then there exists a
    cyclic triangle $x_1\to x_2\to x_3\to x_1$ with $x_i\in X_i$.
\end{corollary}

\begin{corollary}
\label{cor: arcs-large-sets}
    Let $A, B \subseteq \fF_q$.
    If 
    \begin{equation}
    \label{eq: arcs-large-sets-hyp}
        |A||B| \geq 5000q^5,
    \end{equation}
    then there exist flags $a \in A$ and $b \in B$ with $a \to b$ in $D_q$.
\end{corollary}

\noindent To prove the latter corollary, just take $X = A$, $Y = B$, and $Z = \fF_q$ in Proposition~\ref{prp: cyclic-triangles}.

Now let $V$ be a linearly ordered vertex set.
We will define a red-blue coloring of $V^{(3)}$ that depends on an assignment $\phi: V^{(2)} \longrightarrow \fF_q$ of flags to pairs of vertices.
Ultimately, we will take $\phi$ as a uniform random choice from the set of all such maps, but first we define the coloring and show that deterministically, there will not be any red copies of $K_4^-$.

More precisely, color a triple $u < v < w \in V$ red if
\begin{equation}
    \phi(uv) \to \phi(vw) \to \phi(uw) \to \phi(uv),
\end{equation}
and color it blue otherwise.
Notice that this requires not only that the three flags assigned to the three pairs form a cyclic triangle, but that they do so with one of the two possible orientations.

\begin{proposition}
\label{prp: k4-}
    There is no red $K_4^-$ in $V^{(3)}$.
\end{proposition}

\begin{proof}
    If there existed a red $K_4^-$ rooted at a vertex $v \in V$ and containing three other vertices $u_1 < u_2 < u_3 \in V$, a simple case analysis depending on where $v$ appears in the ordering shows that the three red edges among these four vertices always result in a transitive triangle in $D_q$.
    Indeed, if, e.g., $v < u_1 < u_2 < u_3$, the red edges imply
    \[
        \phi(vu_3) \to \phi(vu_2) \to \phi(vu_1), \qquad \phi(vu_3) \to \phi(vu_1).
    \]
    In the remaining three cases $u_1 < v < u_2 < u_3$, $u_1 < u_2 < v < u_3$, and $u_1 < u_2 < u_3 < v$, the corresponding transitive triangles are, respectively,
    \[
    \phi(vu_1) \to \phi(vu_3) \to \phi(vu_2), 
    \qquad 
    \phi(vu_2) \to \phi(vu_1) \to \phi(vu_3), 
    \qquad
    \phi(vu_3) \to \phi(vu_2) \to \phi(vu_1),
    \]
    with an arc from the first to the last flag in each. Each case contradicts Proposition~\ref{prp: Dq-no-TT}.
\end{proof}

Now we would like to study the maximum size of a clique among the blue edges.
For that, we'll use the following consequence of the Saxton--Thomason hypergraph container lemma (specialize Theorem~\ref{thm: k-unif-containers} below to uniformity $h=3$).

\begin{lemma}
\label{lem: 3-unif-containers}
    Let {$0<\e\le 1/288$} be given, and fix a parameter $0 < \tau < 1/2$.
    Suppose $\scH$ is a $3$-uniform hypergraph, $0<e_0\leq e(\scH)/2$, and $\scH$ has the property that every induced subgraph spanning at least $e_0$ edges has average degree $d$ and maximum co-degrees $\Delta_2,\Delta_3$ satisfying
    \begin{equation}
    \label{eq: container-hyp-3}
        \frac{\Delta_2}{d \tau} + \frac{\Delta_3}{d \tau^2} \leq \e.
    \end{equation}
    Then there exists a constant $C = C(\e)$ such that the independent sets of $\scH$ are covered by $\exp \left( C v(\scH) \tau \log (1/\tau) \log (\frac{e(\scH)}{e_0}) \right)$ containers, each spanning at most $e_0$ edges.
\end{lemma}

Fix vertex sets $V_1 < V_2 < V_3 \subseteq V$, each of size $m$, and call an assignment {$\phi: \bigcup_{1\le i<j\le3}(V_i\times V_j) \longrightarrow \fF_q$} {\it triangle-free} (with respect to $V_1 < V_2 < V_3)$ if there do not exist three vertices $v_i \in V_i$, $i \in [3]$, with
\[
    \phi(v_1v_2) \to \phi(v_2v_3) \to \phi(v_1v_3) \to \phi(v_1v_2).
\]
In other words, an assignment is triangle-free if it does not create any transversal red edges between the $V_i$'s.
The following lemma gives a strong bound on the number of triangle-free assignments with respect to a fixed triple of $V_i$'s.

\begin{lemma}
\label{lem: count-three}
    Let $m=M^4$, and suppose $q$ is sufficiently large.
    If $V_1 < V_2 < V_3 \subseteq V$ each have size $m$, then
    the number of triangle-free assignments is at most
    \[
        M^{3m^2-m^2/28}.
    \]
\end{lemma}

\begin{proof}
    Define an auxiliary 3-graph $\scH$ with vertex set $\left( \bigcup_{i < j \in [3]} V_i \times V_j \right) \times \fF_q$ and with, for each $(v_1, v_2, v_3) \in V_1 \times V_2 \times V_3$ and $(c_{12}, c_{23}, c_{13}) \in \fF_q^3$ with $c_{12} \to c_{23} \to c_{13} \to c_{12}$, the triple
    \[
        {((v_1v_2, c_{12}), (v_2v_3, c_{23}), (v_1v_3, c_{13}))} \in E(\scH).
    \]
    In other words, $\scH$ identifies locally the assignments of flags to transversal pairs of vertices that would form a red edge, and the graph of any triangle-free assignment is an independent set of $\scH$.
    We thus want to use the container lemma to bound the number of triangle-free assignments.
    Observe that
    \[
        \Delta_2(\scH) \leq M \qquad \Delta_3(\scH) \leq 1 \qquad v(\scH) = 3 m^2 M \qquad e(\scH) \leq m^3 M^3,
    \]
    
    and set
    \[
        e_0:=\frac{m^3}{4},\qquad
        \tau:=K\sqrt{M/m}=KM^{-3/2},
    \]
    where $K$ is a sufficiently large constant. We apply Lemma~\ref{lem: 3-unif-containers} with $\e=1/288$.
    By Proposition~\ref{prp: cyclic-triangles}, every transversal
    triple of vertices gives rise to at least one edge of $\scH$
    for sufficiently large $q$, so $e(\scH)\geq m^3>e_0$.
    For any $U\subseteq V(\scH)$ with $e(\scH[U])\geq e_0$,
    the average degree of $G=\scH[U]$ satisfies
    \[
        d(G)\geq\frac{3e_0}{3m^2M}=\frac{m}{4M},
    \]
    and thus
    \begin{align*}
        \frac{\Delta_2(G)}{d(G)\tau}
              +\frac{\Delta_3(G)}{d(G)\tau^2}
        &\leq \frac{4M^{7/2}}{Km}+\frac{4M^4}{K^2m}\\
        &=\frac{4}{KM^{1/2}}+\frac{4}{K^2}.
    \end{align*}
    Choosing $K$ sufficiently large and then $q$ sufficiently
    large ensures that this is less than $\e$,
    verifying~\eqref{eq: container-hyp-3}.

    Thus, Lemma~\ref{lem: 3-unif-containers} applies, and
    the independent sets of $\scH$ are covered by a family
    $\scC$ of containers, each spanning at most $e_0$ edges.
    Since $v(\scH)=3m^2M$, $e(\scH)\leq m^3M^3$, and
    $m=M^4$, the number of containers satisfies
    \begin{align*}
        \log|\scC|
        &\leq C_3v(\scH)\tau\log(1/\tau)
                      \log\left(\frac{e(\scH)}{e_0}\right)\\
        &=O\bigl(M^{15/2}(\log M)^2\bigr)
         =o(m^2\log M).
    \end{align*}
    In particular, for sufficiently large $q$,
    $|\scC|\leq M^{m^2/28}$.

    Now let $C$ be one such container.
    Call a transversal pair $e \in V_i \times V_j$, $i \neq j$, {\it small} if
    \[
        |\{ c \in \fF_q: (e,c) \in C \}| < M^{6/7},
    \]
    and call it {\it large} otherwise.
    Note that if less than $m^2/2$ pairs are small, then less than $m \cdot m^2/2$ transversal triples of vertices contain a small pair, and so at least $m^3/2$ of them have all pairs large;
    but whenever this happens, by Corollary~\ref{cor: cyclic-triangles}, this triple has an assignment in $C$ which yields an edge of $\scH$, and so we would have $e(\scH[C]) \geq m^3/2 > e_0$, which is impossible.

    Therefore, at least $m^2/2$ pairs are small.
    But then $C$ can only support at most $M^{3m^2-m^2/14}$ triangle-free assignments.
    Taking a union over all containers $C$ and using the bound obtained above on their number gives the desired result.
\end{proof}

Now we are ready to prove the claimed bound on $r_3(K_4^-, K_n)$.

\begin{corollary}
    For every sufficiently large prime power $q$, with $m=M^4$ and $n=3m$, we have
    \[
        r_3(K_4^-,K_n)\geq M^{m/168}.
    \]
\end{corollary}

\begin{proof}
    Let $N=M^{m/168}$, and take
    $\phi:[N]^{(2)}\longrightarrow\fF_q$ uniformly at random.
    Color each triple red or blue based on the rule above.
    By Proposition~\ref{prp: k4-}, there is (deterministically)
    no red $K_4^-$ under this coloring.
    Suppose, then, that $W$ is a set of $n=3m$ vertices,
    and let $V_1<V_2<V_3$ be the first, middle, and last
    $m$ vertices of $W$. Then Lemma~\ref{lem: count-three}
    gives
    \begin{align*}
        \bbP[W\text{ is an all-blue clique}]
        &\leq\bbP[V_1\times V_2\times V_3\text{ is all-blue}]\\
        &\leq\frac{M^{3m^2-m^2/28}}{M^{3m^2}}
         =M^{-m^2/28}.
    \end{align*}
    Since there are at most $N^{3m}$ choices for $W$,
    a union bound gives that the probability that there exists
    an all-blue clique of size $n$ is at most
    \[
        N^{3m}M^{-m^2/28}
        \leq M^{m^2/56-m^2/28}<1.
    \]
    Thus, there exists a choice of $\phi$ which produces
    a red-blue coloring with no red $K_4^-$ and no blue $K_n$.
\end{proof}

To obtain the bound for every sufficiently large $n$, choose
an odd prime $q=\Theta(n^{1/12})$ with the implicit constant
sufficiently small that $3M^4\leq n$. This is possible by
Bertrand's postulate, since $M=\Theta(q^3)$.
With $m=M^4$, we then have $3m\leq n$ and $m=\Theta(n)$.
The bound just obtained for $3m$, together with monotonicity
in the blue clique size, gives
\[
    r_3(K_4^-,K_n)\geq M^{m/168}
    =2^{\Omega(n\log n)}.
\]


\section{A general framework for bounding the size of the largest blue clique}
\label{sec: ind-number}

In this section, we focus on the common probabilistic ingredient in the proofs of
Theorems~\ref{thm: 3-edges-k-unif} and~\ref{thm: half-graph}, which we believe
could be of independent interest. In each construction, we independently assign
labels from a finite set $\Omega$ to certain $d$-subsets of an ordered vertex
set, and color a $k$-edge red whenever its labels form one of a prescribed
collection of patterns. To bound the size of the largest blue clique, it thus
suffices to control the probability that no red edge appears between $k$ ordered vertex sets.

We will show that the following local condition on the patterns is enough for such a bound:
if every Cartesian product of sufficiently large subsets of $\Omega$ contains
one of the prescribed patterns, then the number of assignments avoiding all such patterns
between the $k$ ordered vertex sets is exponentially smaller than the total
number of assignments. 
We prove this using the hypergraph container method
of Saxton and Thomason. 
One particular advantage of this formulation is that the probabilistic argument is independent of the particular forbidden hypergraph, allowing us to use the same counting lemma in both constructions of Section~\ref{sec: proofs-k-unif}.

More formally, let $k \geq 3$ and $d < k$ be positive integers, and let $\scF \subseteq [k]^{(d)}$ be a family of $d$-subsets of $[k]$ which cover $[k]$.
We will work with a totally ordered vertex set $V$ along with a ``target set'' $\Omega$ of size $T$.

Suppose that we have a map
$\varphi: V^{(d)}\longrightarrow\Omega$ and a nonempty collection
$\scP\subseteq\Omega^{\scF}$ of ``patterns.''

For a given $k$-subset $e = \{v_1 < \dots < v_k\}$ and $F \in \scF$, write $e_F = \{v_i: i \in F\}$.
We will then define a $k$-uniform hypergraph $(V, E)$ on vertex set $V$ with
\[
    e \in E \Longleftrightarrow \langle \varphi(e_F) \rangle_{F \in \scF} \in \scP.
\]
Our goal will be to understand the independence number of this hypergraph given some information about $\varphi$ and $\scP$.

\subsection{Containers}

We will eventually use a version of the Saxton--Thomason hypergraph container theorem~\cite{ST15}, but to do so we need to review some of the concepts involved in their result.
To that end, let $\scH$ be an $h$-uniform hypergraph.
The usual degree of a subset $J \subseteq V(\scH)$ will be denoted by $d(J) = |\{e \in E(\scH): J \subseteq e\}|$, and for a vertex $v \in V(\scH)$, we will set its maximum $j$-degree to be
\[
    d_j(v) := \max \left\{ d(J): v \in J, |J|=j \right\}.
\]
The Saxton--Thomason co-degree function is then defined (with respect to a formal variable $\tau$) by
\begin{equation}
\label{eq: ST-codegree}
    \delta(\scH,\tau)
    =
    \frac{2^{\binom{h}{2}-1}}{v(\scH) \cdot d(\scH)} \sum_{j=2}^h \frac{\sum_{v \in V(\scH)} d_j(v)}{2^{\binom{j-1}{2}} \tau^{j-1}},
\end{equation}
where 
$d(\scH) = \frac{1}{v(\scH)} \sum_{v \in V(\scH)} d(v)$ denotes the average vertex-degree.

The following container theorem is an abridged version of~\cite[Theorem 6.3]{ST15} due to Saxton and Thomason.
We have specialized their result to constant $\tau$, suppressed parts of the conclusion which we do not use (particularly the parts about the $T$'s in (a) and the parts after ``Moreover"), and uniformly applied the bound $v(\scH) \, \tau \log( 1/\tau)$ to each of the terms in the sum in part (c), the number of which is $O_h(\log(e(\scH)/e_0))$.

\begin{theorem}
[Saxton--Thomason]
\label{thm: k-unif-containers}
    Let $\scH$ be a $h$-uniform hypergraph, and let $0<e_0 \leq e(\scH)/2$.
    For a given $0 < \tau < 1/2$, suppose that the co-degree function satisfies
    \begin{equation}
    \label{eq: container-hypothesis}
        \d(\scH[U], \tau) \leq \frac{1}{12h!}
    \end{equation}
    for all vertex subsets $U \subseteq V(\scH)$ with $e(\scH[U]) \geq e_0$.
    Then there is a family $\scC$ of subsets of $V(\scH)$ such that the following holds:
    \begin{enumerate}
        \item every independent set of $\scH$ is contained in some $C \in \scC$;
        \item for all $C \in \scC$, we have $e(\scH[C]) 
        \leq e_0$;
        \item there is a constant $C_h > 0$, depending only on $h$, such that
        \[
        |\scC|
        \leq \exp \left( C_h \, v(\scH) \, \tau \log(1/\tau) \log \left( \frac{e(\scH)}{e_0} \right) \right).
        \]
    \end{enumerate}
\end{theorem}

As in the proof that $r_3(K_4^-, K_n) \geq \exp(\Omega(n \log n))$, our goal will be to apply this to an auxiliary hypergraph below.

\subsection{An assignment-counting lemma}
\label{section: key-lemma}

Now let $V_1 < V_2 < \cdots < V_k \subseteq V$ be disjoint subsets of our ordered vertex set $V$, each of size $m$.
Put
\[
    \scX:=\bigcup_{F \in \scF} \prod_{i \in F} V_i.
\]
In other words, $\scX$ is the set of all transversal $d$-sets of vertices from the $V_i$'s whose indices correspond to some $F \in \scF \subset [k]^{(d)}$; note that 
\[
  |\scX|=|\scF| \cdot m^d.
\]
For a transversal $\mathbf{v} = (v_1, \dots, v_k) \in V_1 \times \dots \times V_k$, we will write $\mathbf{v}_F = (v_i)_{i \in F}$ to denote a subsequence of it which lies in $\scX$.
Finally, given a set of patterns $\scP \subseteq \Omega^\scF$ as above, we will say that a map 
$\varphi: \scX \longrightarrow \Omega$
is a {\it $\scP$-free assignment} if for all $\mathbf{v} \in \prod_{i=1}^k V_i$, we have
\[
  (\varphi(\mathbf{v}_F))_{F \in \scF} \notin \scP.
\]
Following the heuristic given at the beginning of the section, the set $V_1 \cup \dots \cup V_k$ is {independent in the hypergraph on $V$ only if} $\varphi: \scX \longrightarrow \Omega$ is $\scP$-free.
Our key lemma below bounds the number of such $\scP$-free assignments, and shows in a precise sense that ``most'' maps $\varphi$ do not result in $V_1 \cup \dots \cup V_k$ being independent in the hypergraph on $V$.


\begin{lemma}
\label{lem: assignment-counting}
    Fix integers $m \geq 1$, $k \geq 3$, and $1 \leq d < k$, and let $\scF \subseteq [k]^{(d)}$ be a family of $|\scF| = h \geq 2$ sets, the union of which covers $[k]$.
    Fix also constants $0 < \eta < 1$ and
    \[
        0 < b \leq a < \frac{1}{h}.
    \]
    Suppose that for a set $\Omega$ of size $T$, where
    \[
        m^b \leq T \leq m^a,
    \]
    and for a set $\scP \subseteq \Omega^\scF$ of patterns,
    we have the property:
    \begin{itemize}[label = $(\star)$]
        \item \label{it: box-property} every $\prod_{F \in \scF} W_F \subseteq \Omega^{\scF}$ with $|W_F| \geq T^{1-\eta}$ for each $F \in \scF$ must
        contain an element of $\scP$.
    \end{itemize}
    Finally, let $V_1, \dots, V_k$ be {pairwise disjoint, ordered} sets of size $m$ and set
    \[
        \scX = \bigcup_{F \in \scF} \prod_{i \in F} V_i.
    \]
    Then if $m$ is sufficiently large, the number of $\scP$-free assignments $\varphi: \scX \longrightarrow \Omega$ is at most
    \[
        T^{hm^d-(\eta/4)m^d}.
    \]
\end{lemma}

The basic idea of the proof is to construct an auxiliary hypergraph with edges identifying potential local assignments that would create a transversal $k$-edge with an assignment in $\scP$, thus putting $\scP$-free assignments into correspondence with independent sets in this auxiliary hypergraph.
The hypothesis~$(\star)$ and the resulting symmetry in the auxiliary hypergraph are all we will need to invoke the container lemma, which will give the claimed strong bound on the number of $\scP$-free assignments.

\begin{proof}
    We define an auxiliary $h$-uniform hypergraph $\scH$ with vertex set
    \[
        V(\scH) = 
        \scX \times \Omega.
    \]
    For every $\mathbf{v} \in \prod_{i=1}^k V_i$ and
    $(\omega_F)_{F \in \scF} \in \scP$,
    put in $\scH$ the hyperedge
    \[
        \{(\mathbf{v}_F, \omega_F): F \in \scF\}.
    \]
    An assignment 
    $\varphi : \scX \longrightarrow \Omega$ 
    is then $\scP$-free precisely
    when
    \[
        I_\varphi := \{(\mathbf{x}, \varphi(\mathbf{x})): \mathbf{x} \in \scX\}
    \]
    is an independent set of $\scH$, and so to bound the number of $\scP$-free assignments it suffices to count the independent sets of $\scH$ that arise in this way.

    This auxiliary hypergraph $\scH$ has
    \[
        v(\scH) = h m^d \cdot T
        \qquad
        \text{and}
        \qquad
        m^k \leq e(\scH) \leq m^k \cdot T^h
    \]
    (since $\scF$ covers $[k]$, distinct $\mathbf{v}, \mathbf{v}' \in \prod_{i=1}^k V_i$ must give distinct sequences $(\mathbf{v}_F)_{F \in \scF}, (\mathbf{v}'_F)_{F \in \scF}$),
    and for every {$2 \leq j \leq h$},
    \begin{equation}
    \label{eq: k-unif-codeg}
        \Delta_j(\scH) \leq m^{k-d-1} T^{h-j};
    \end{equation}
    indeed, if $j \geq 2$ auxiliary vertices lie in a common edge, then at least $d+1$ positions are determined in the
    underlying transversal through the $V_i$'s, and so it only remains to choose the remaining $k-d-1$ positions in this transversal along with the corresponding $\omega_F$'s from $\Omega$.

    Set
    \[
        e_0 := \frac{m^k}{4}
        \qquad
        \text{and}
        \qquad
        \tau := K \left(\frac{T}{m}\right)^{1/(h-1)}
    \]
    where $K$ is a sufficiently large constant which will ultimately only depend on $h$;
    note that since $T \leq m^a$ and $a < 1/h$, we have 
    $\tau < 1/2$ for sufficiently large $m$.

    Now we would like to verify that the hypothesis~\eqref{eq: container-hypothesis} holds for this choice of $e_0$ and $\tau$ in $\scH$.
    Let $U \subseteq V(\scH)$ satisfy $e(\scH[U]) \geq e_0$, and let $d_U$ be the average degree of $\scH[U]$.  
    Then
    \[
        d_U
        = \frac{h \cdot e(\scH[U])}{|U|}
        \geq \frac{h m^k/4}{h m^d \cdot T}
        =\frac{m^{k-d}}{4T}.
    \]
    For $2 \leq j \leq h$, \eqref{eq: k-unif-codeg} gives
    \[
        \sum_{u\in U} d_j(u)
        \leq |U| m^{k-d-1} T^{h-j}.
    \]
    Consequently,
    \begin{align*}
        \d (\scH[U],\tau)
        &\leq \frac{2^{\binom{h}{2}-1}}{|U| \cdot d_U} \sum_{j=2}^h \frac{\sum_{u \in U} d_j(u)}{2^{\binom{j-1}{2}} \tau^{j-1}}\\
        &\leq \frac{2^{\binom{h}{2}-1}}{|U| \cdot \frac{m^{k-d}}{4T}} \sum_{j=2}^h \frac{|U| m^{k-d-1} T^{h-j}}{\tau^{j-1}}\\
        &= 2^{\binom{h}{2}-1} \sum_{j=2}^h \frac{4T}{m} \cdot T^{h-j} \cdot K^{-(j-1)} \left( \frac{m}{T} \right)^{\frac{j-1}{h-1}}\\
        &= 4 \cdot 2^{\binom{h}{2}-1} \sum_{j=2}^h T^{h-j+1 - \frac{j-1}{h-1}} m^{\frac{j-1}{h-1}-1} \cdot K^{-(j-1)}\\
        &= 4 \cdot 2^{\binom{h}{2}-1} \sum_{j=2}^h \left( \frac{T^h}{m} \right)^{\frac{h-j}{h-1}} \cdot K^{-(j-1)}
    \end{align*}
    Since $T \leq m^a$ and $a < 1/h$, the $j$th term in the sum is at most $K^{-(j-1)}$; choosing $K$ sufficiently large with respect to $h$ thus gives
    \[
        \d (\scH[U],\tau) \leq \frac{1}{12h!}
    \]
    for every such $U$, as required.

    Thus, Theorem~\ref{thm: k-unif-containers} applies.  
    The independent sets of $\scH$ are covered by a family $\scC$ satisfying
    \[
        e(\scH[C]) \leq e_0
    \]
    for all $C \in \scC$ and
    \begin{align*}
        |\scC|
        &\leq \exp \left( C_h \, v(\scH) \, \tau \log(1/\tau) \log \left( \frac{e(\scH)}{e_0} \right) \right)\\
        &\leq \exp \left( C_h \cdot h m^d T \cdot \left( \frac{T}{m} \right)^{1/(h-1)} \log(m/T) \cdot \log \left( \frac{m^k T^h}{m^k/4} \right) \right)\\
        &\leq \exp \left( C_h \cdot h m^{d - \frac{1}{h-1}} T^{\frac{h}{h-1}} \log(m^{1-b}) \cdot \log (4 T^h) \right)\\
        &\leq \exp \left( C_{h, b} \cdot m^{d - \frac{1}{h-1}} m^{a \cdot \frac{h}{h-1}} \log m \cdot \log T \right).
    \end{align*}
    Since $a < 1/h$, it follows that as $m$ gets large,
    \[
        |\scC| 
        \leq \exp \left( o \left( m^d \log T \right) \right).
    \]

    Now fix a container $C \in \scC$. 
    For every $\mathbf{x} \in \scX$, define its {\it list} to be the set
    \[
        L_{\mathbf{x}} := \{\omega \in \Omega: (\mathbf{x}, \omega) \in C\}.
    \]
    Call $\mathbf{x}$ {\it small} if its list has size
    \[
        |L_{\mathbf{x}}| < T^{1-\eta},
    \]
    and {\it large} otherwise.
    
    \begin{claim}
        More than $m^d/2$ elements $\mathbf{x} \in \scX$ are small.
    \end{claim}
    
    \begin{proof}[Proof of Claim]
        Assume otherwise.
        Then since every $d$-tuple $\mathbf{x} \in \scX$ belongs to exactly $m^{k-d}$ transversal $k$-tuples in $\prod_{i=1}^k V_i$, at most
        \[
            \frac{m^d}2\cdot m^{k-d} =\frac{m^k}{2}
        \]
        transversal $k$-tuples contain a small $\mathbf{x}$. Hence, at least $m^k/2$ transversal $k$-tuples $\mathbf{v}$ satisfy
        \[
            |L_{\mathbf{v}_F}|
            \geq T^{1-\eta}
            \qquad\text{for every } F \in \scF.
        \]
        For each such $\mathbf{v}$, the hypothesis~$(\star)$ supplies a member of
        \[
            \scP \cap
            \prod_{F \in \scF} L_{\mathbf{v}_F},
        \]
        and hence an edge of $\scH[C]$.  
        Therefore
        \[
            e(\scH[C]) \geq \frac{m^k}{2} > e_0,
        \]
        which is a contradiction.
    \end{proof} 
    
    Thus, more than $m^d/2$ lists are small.  
    The number of assignments $\varphi: \scX \longrightarrow \Omega$ supported by $C$ is then at most
    \[
        \prod_{\mathbf{x} \in \scX} |L_{\mathbf{x}}|
        \leq
        T^{h m^d - (\eta/2)m^d}.
    \]
    Since every $\scP$-free assignment corresponds to a unique independent set in $\scH$, taking $m$ large enough so that we can bound
    \[
        |\scC| \leq \exp \left( (\eta / 4) m^d \log T \right)
    \]
    then gives the claimed bound on the number of $\scP$-free assignments.
\end{proof}


\section{Forbidding fixed hypergraphs and the proofs of Theorems~\ref{thm: 3-edges-k-unif} and~\ref{thm: half-graph}}
\label{sec: proofs-k-unif}

In view of Lemma~\ref{lem: assignment-counting}, the main task is now to construct collections of patterns satisfying the
Cartesian-product condition~$(\star)$. In both applications, our label set
will be $\Omega=\fF_q^d$ for a suitable choice of $d$, where $\fF_q$ is the
vertex set of the directed graph $D_q$ introduced in Section~\ref{sec: K4-}.
The patterns will be specified by requiring certain coordinates of these
labels to form directed edges in $D_q$.

We will show that this condition is guaranteed whenever the directed
constraints form a vertex-disjoint union of out-stars. The reason is that
Corollary~\ref{cor: arcs-large-sets} allows us to satisfy the constraints
one star at a time.
This yields a general
lemma which we will apply twice: first when the directed
constraints form a matching (for Theorem~\ref{thm: 3-edges-k-unif}), and then when they form a union of larger
out-stars (for Theorem~\ref{thm: half-graph}).

Given positive integers $r$ and $d$, we will call a directed graph $G$ on vertex set $[r] \times [d]$ a {\it pattern digraph} if:
\begin{itemize}
    \item $G$ is a vertex-disjoint union of directed out-stars (single directed edges are allowed); and
    \item each such star contains at most one vertex from each row $\{i\} \times [d]$, $i \in [r]$.
\end{itemize}

\begin{lemma}
\label{lem: large-box}
    Fix positive integers $r$ and $d$, and let $G$ be a pattern digraph on $[r] \times [d]$.

    For every $0 < \eta < 1/6$ and all sufficiently large prime powers $q$, we have the following: if $X_1, \dots, X_r \subseteq \fF_q^d$ have $|X_\l| \geq M^{d - \eta}$ for each $\l \in [r]$, then there exist
    \[
        x_i = (x_i^{(1)}, \dots, x_i^{(d)}) \in X_i
    \]
    such that 
    \[
        x_i^{(a)} \longrightarrow x_j^{(b)} \text{ in } D_q \qquad \text{ whenever } \qquad (i,a) \longrightarrow (j,b) \text{ in } G.
    \]
\end{lemma}

\begin{proof}
    We greedily find each $x_i \in X_i$, making repeated use of Corollary~\ref{cor: arcs-large-sets}.

    \begin{claim*}
    \label{clm: out-star}
        For all $c > 0$ and $s \in \N$, the following holds for $q$ sufficiently large.
        Given $A, B_1, \dots, B_s \subseteq \fF_q$ which all have size at least $cM^{1-\eta}$, we can always find a directed out-star $a \to b_i$ with $a \in A$ and $b_i \in B_i$ for $1 \leq i \leq s$.
    \end{claim*}

    \begin{proof}[Proof of Claim~\ref{clm: out-star}]
        By Corollary~\ref{cor: arcs-large-sets}, at most $5000q^5/|B_i|$ flags in $A$ have no out-neighbors in $B_i$, and so at most
        \[
            \sum_{i=1}^s \frac{5000 q^5}{|B_i|} \leq s \cdot \frac{5000 q^5}{cM^{1-\eta}}
        \]
        flags in $A$ fail to have an out-neighbor in every $B_i$; since $\eta < 1/6$ and $M \geq q^3$, this quantity will be less than $cM^{1-\eta} \leq |A|$ for $q$ sufficiently large, and so there must exist an $a \in A$ with an out-neighbor in every $B_i$.
    \end{proof}

    We greedily construct the $x_i \in X_i$ one star in $G$ at a time using Claim~\ref{clm: out-star}.
    The order in which each star is addressed is insignificant, and thus may be done arbitrarily, as the stars are vertex-disjoint.
    Suppose we are at the step addressing a star in $G$ of the form $(i, a) \to (j_1, b_1), \dots, (i,a) \to (j_s, b_s)$.
    For $\l \in [r]$, let
    \begin{itemize}
        \item $S_\l \subseteq [d]$ be the set of coordinates which have already been fixed in $x_\l$, and $s_\l = |S_\l|$;
        \item $D_\l = [d] \setminus S_\l$ be the set of coordinates remaining to be chosen in $x_\l$, and $d_\l = |D_\l|$;
        \item $Y_\l \subseteq \fF_q^{D_\l}$ be the surviving set of $d_\l$-tuples of flags that may be inserted into the remaining coordinates $D_\l \subseteq [d]$ of $x_\l$ to obtain an element of $X_\l \subseteq \fF_q^{S_\l \dot\cup D_\l}$; and
        \item $Y_\l^{(k)} = \{x^{(k)}: x \in Y_\l\}$ be the flags which appear as the $k$th coordinate of some vector in $Y_\l$, where $k \in D_\l$.
    \end{itemize}
    The inductive hypothesis we will maintain is that 
    \begin{equation}
    \label{eq: Y_l-size}
        |Y_\l| \geq 2^{-s_\l} M^{d_\l - \eta}
    \end{equation}
    for all $\l \in [r]$.
    (Note that this follows from the assumption that each $|X_\l| \geq M^{d-\eta}$ at the initial step.)
    Also observe that if we take, e.g., $q \geq 2^{\frac{d}{3(1-\eta)}}$, then $M \geq q^3$ guarantees that the sets $Y_\l$ will always be nonempty when $d_\l > 0$.

    Now, since we are at the step corresponding to the star $(i, a) \to (j_1, b_1), \dots, (i,a) \to (j_s, b_s)$ in $G$, we would like to choose $x_i^{(a)} \in Y_i^{(a)}$ and $x_{j_1}^{(b_1)} \in Y_{j_1}^{(b_1)}, \dots, x_{j_s}^{(b_s)} \in Y_{j_s}^{(b_s)}$ with $x_i^{(a)} \to x_{j_1}^{(b_1)}, \dots, x_i^{(a)} \to x_{j_s}^{(b_s)}$ and which preserve the hypothesis~\eqref{eq: Y_l-size}.
    To that end, for $z \in \fF_q$, let $Y_i^{(a)}(z)$ denote the updated set {$Y_i$} if $z$ were inserted into $x_i^{(a)}$, that is, the set of vectors in $\fF_q^{D_i \setminus \{a\}}$ which, when inserted into positions $D_i \setminus \{a\}$ and $z$ inserted into position $a$, result in an element of $X_i$.
    The sets $Y_{j_1}^{(b_1)}(z), \dots, Y_{j_s}^{(b_s)}(z)$ are defined analogously.
    We will call $z$ {\it heavy with respect to $Y_i^{(a)}$} if 
    \begin{equation}
    \label{eq: heavy}
        |Y_i^{(a)}(z)| \geq |Y_i|/2M,
    \end{equation}
    and {\it heavy with respect to $Y_{j_k}^{(b_k)}$} if the analogous inequality holds for $Y_{j_k}^{(b_k)}(z)$.
    Finally, let $H_i^{(a)}$ and $H_{j_1}^{(b_1)}, \dots, H_{j_s}^{(b_s)}$ denote the set of heavy flags with respect to $Y_i^{(a)}$ and $Y_{j_1}^{(b_1)}, \dots, Y_{j_s}^{(b_s)}$, respectively.

    Crucially, observe that $|H_i^{(a)}|, |H_{j_1}^{(b_1)}|, \dots, |H_{j_s}^{(b_s)}| \geq 2^{-d} M^{1-\eta}$: indeed, if $z \notin H_i^{(a)}$, for instance, then it accounts for at most $|Y_i|/2M$ elements of $Y_i$, while if $z \in H_i^{(a)}$, it accounts for at most $M^{d_i-1}$ such elements; using that there are $M$ flags in total, this gives
    \[
        |Y_i| \leq |H_i^{(a)}| \cdot M^{d_i-1} + M \cdot \frac{|Y_i|}{2M},
    \]
    and rearranging this and applying~\eqref{eq: Y_l-size} gives the claimed bound.
    
    Invoking Claim~\ref{clm: out-star} with $c = 2^{-d}$, $A = H_i^{(a)}$, and $B_k = H_{j_k}^{(b_k)}$ for $k \in [s]$, it follows that for $q$ sufficiently large there must exist $x_i^{(a)} \in H_i^{(a)} \subseteq Y_i^{(a)}$ and $x_{j_k}^{(b_k)} \in H_{j_k}^{(b_k)} \subseteq Y_{j_k}^{(b_k)}$ with $x_i^{(a)} \to x_{j_k}^{(b_k)}$ in $D_q$ for all $k \in [s]$; and moreover, by choosing heavy flags from all these sets,~\eqref{eq: heavy} guarantees that the inductive hypothesis~\eqref{eq: Y_l-size} is preserved. Since each star meets each row at most once, at most one coordinate of any $x_\l$ is fixed at this step.
    After processing all stars, choose any surviving completion in each $X_\l$; the positive lower bound in~\eqref{eq: Y_l-size} guarantees that such a completion exists.

\end{proof}

\subsection{Forbidding three red \texorpdfstring{$k$-edges among $k+1$}{k-edges among k+1} vertices}

Let $V$ be a totally ordered vertex set.
Our goal will be to define a red-blue coloring of $V^{(k)}$ so that no $(k+1)$-set of vertices spans three red edges, but the size of the largest blue clique is still controlled.

To that end, let $q$ be a prime power (which we will choose later) and fix a map $\varphi: V^{(k-1)} \longrightarrow \fF_q^{k-1}$. 
(We will ultimately pick $\varphi$ uniformly at random among all such maps.)
For a $(k-1)$-set $f = \{v_1 < \dots < v_{k-1}\}$ and $i \in [k-1]$, we will write
\[
    \varphi(f; v_i) := (\varphi(f))_i;
\]
that is, if $v$ is the $i$th vertex of the $(k-1)$-set $f$ under the ordering $<$, then $\varphi(f; v)$ is the $i$th flag in $\varphi(f) \in \fF_q^{k-1}$.
We color a $k$-edge $e = \{v_1 < \cdots < v_k\}$
red if
\begin{equation}
\label{eq: red-edges}
    \varphi(e \setminus \{v_j\}; v_i) \longrightarrow \varphi(e \setminus \{v_i\}; v_j) \text{ in } D_q
    \qquad
    \text{for every } i < j \in [k];
\end{equation}
otherwise, we color it blue.

\begin{lemma}
\label{lem: forbidden-red-config}
    For every assignment
    \[
    \varphi: V^{(k-1)} \longrightarrow {\fF_q^{k-1}},
    \]
    every $(k+1)$-subset of $V$ spans at most two red edges.
\end{lemma}

\begin{proof}
    Suppose, to the contrary, that a $(k+1)$-set $F$ spans three red edges, and let $u < v < w$ be three vertices {such that removing any one of them results in a red edge}.
    In other words, if $S = F \setminus \{u, v, w\}$, then the edges $S \cup \{u,v\}$, $S \cup \{v, w\}$, and $S \cup \{u, w\}$ are all red.
    Consider the three flags $\varphi(S \cup \{u\}; u)$, $ \varphi(S \cup \{v\}; v)$, and $\varphi(S \cup \{w\}; w)$.
    By applying \eqref{eq: red-edges} to the red edge $S \cup \{x, y\}$ for each $x < y \in \{u, v, w\}$, we get that
    \[
        \varphi(S \cup \{x\}; x) \to \varphi(S \cup \{y\}; y)
    \]
    for each such pair of flags.
    But this means they form a transitive triangle in $D_q$, which is a contradiction to Proposition~\ref{prp: Dq-no-TT}.
\end{proof}

Now consider a $k$-tuple of $(k-1)$-tuples of flags $(x_1, \dots, x_k) \in (\fF_q^{k-1})^k$, and write\\
$x_i =: (x_i^{(1)}, \dots, x_i^{(i-1)}, x_i^{(i+1)}, \dots, x_i^{(k)})$.
That is, we label the flags in the $k-1$ coordinates of $x_i$ in superscript by the elements of $[k] \setminus \{i\}$ in order.
We will call such a tuple $(x_1, \dots, x_k) \in (\fF_q^{k-1})^k$ a {\it red pattern} if
\begin{equation}
\label{eq: red-pattern}
    x_j^{(i)} \to x_i^{(j)} \text{ in } D_q \qquad \text{ for all } i < j \in [k].
\end{equation}
Observe that, given an assignment $\varphi: V^{(k-1)} \longrightarrow \fF_q^{k-1}$, an edge $e = \{v_1 < \dots < v_k\}$ is colored red if and only if $\langle \varphi(v_1, \dots, \hat{v_i}, \dots, v_k) \rangle_{i=1}^k$ is a red pattern.

\begin{lemma}
\label{lem: 3-edge-large-box}
    For every $0 < \eta < 1/6$ and all sufficiently large prime powers $q$, we have the following: if $X_1, \dots, X_k \subseteq \fF_q^{k-1}$ have $|X_\l| \geq M^{k-1 - \eta}$ for each $\l \in [k]$, then $\prod_{\l=1}^k X_\l$ contains a red pattern.
\end{lemma}

\begin{proof}
    Apply Lemma~\ref{lem: large-box} with $r=k$ and $d=k-1$, and $G$ defined on $[k] \times [k-1]$ as follows: for each row $\{i\} \times [k-1]$, index the elements by $[k] \setminus \{i\}$ following the order of the second coordinate; then take $G$ to be the union of vertex-disjoint directed edges
    \[
        (j,i) \longrightarrow (i,j) \qquad \text{for }1\leq i<j\leq k.
    \]
\end{proof}

We will prove Theorem~\ref{thm: 3-edges-k-unif} together with
Theorem~\ref{thm: half-graph} at the end of the section, using the same
assignment-counting argument for both constructions.


\subsection{Forbidding a half-graph}
As above, our goal will be to define a red-blue coloring of $[N]^{(k)}$ for suitable $N$ with neither a red copy of $B^{(k)}$ nor a blue copy of $K_n^{(k)}$.

Fix a prime power $q$ and an assignment $\varphi: [N]^{(2)} \longrightarrow \fF_q^2$ from pairs of vertices to pairs of flags.
For a pair $\{v_1,v_2\}\in [N]^{(2)}$ with $v_1<v_2$, we will again write
\[
    \varphi(v_1v_2; v_i) = (\varphi(v_1v_2))_i.
\]
With this notation, we will color a $k$-edge $e = \{v_1 < \dots < v_k\}$ red if:
\begin{itemize}
    \item $\varphi(v_1v_2; v_1) \longrightarrow \varphi(v_1v_i; v_1)$ for all $2 < i \leq k$; and
    \item $\varphi(v_1v_i; v_i) \longrightarrow \varphi(v_iv_j; v_i)$ in $D_q$ for all $2\leq i\leq k$ and $j \in [k] \setminus \{1, i\}$,
\end{itemize}
and otherwise we color $e$ blue.
In words, we color $e$ red if, for all $i \in [k]$, when we fix the index $\varphi(-,v_i; v_i)$, the result of plugging in $v_j$ with $j \neq i$ is a directed out-star in $D_q$ rooted at the vertex with the smallest index not equal to $i$.
We will call this out-star the {\it out-star indexed by $v_i$}.

\begin{lemma}
\label{lem: forbidden-halfgraph}
    For any $q$ and $N$ and any assignment $\varphi: [N]^{(2)} \longrightarrow \fF_q^2$, there is no red copy of $B^{(k)}$ under the above rule.
\end{lemma}

\begin{proof}
    Suppose there were a red copy of $B^{(k)}$ on disjoint {$(k-1)$-sets} $S$ and $T$ with distinguished vertex $s_* \in S$.
    Recall that the edges are
    \[
        \{S \cup \{t\}: t \in T\} \,\, \cup \,\, {\{T \cup \{s_*\}\}}.
    \]
    Let $C = S \setminus \{s_*\}$ and set
    \[
        c_0 := \min_< C \qquad t_0 := \min_< T.
    \]
    Observe that we have the following arcs in $D_q$:
    \begin{enumerate}
        \item \label{it: halfgraph-1} since $T \cup \{s_*\}$ is red, $\varphi(s_* t_0; s_*) \to \varphi(s_* t; s_*)$ for all $t \in T \setminus \{t_0\}$;
        \item \label{it: halfgraph-2} since $S \cup \{t\}$ is red for all $t \in T$, we have $\varphi(s_* c_0; s_*) \to \varphi(s_* t; s_*)$ if $c_0 < t$, and $\varphi(s_* t; s_*) \to \varphi(s_* c_0; s_*)$ if $t < c_0$.
    \end{enumerate}
    In other words, choose $t_1\in T\setminus\{t_0\}$, which exists since $k\geq3$. The three flags $\varphi(s_*c_0;s_*)$, $\varphi(s_*t_0;s_*)$, and $\varphi(s_*t_1;s_*)$ form a transitive triangle in $D_q$, since the directions of all three arcs follow the ordering of $c_0,t_0,t_1$. This is a contradiction.
\end{proof}

We will call a $\binom{k}{2}$-tuple of pairs of flags $(x_{1,2}, \dots, x_{k-1,k}) \in (\fF_q^2)^{\binom{k}{2}}$ a {\it half-graph pattern} if:
\begin{itemize}
    \item $(x_{1,2})_1 \longrightarrow (x_{1,i})_1$ for all $2 < i \leq k$;
    \item $(x_{1,i})_2 \longrightarrow (x_{j,i})_2$ in $D_q$ for all $1 < j < i$; and
    \item $(x_{1,i})_2 \longrightarrow (x_{i,j})_1$ in $D_q$ for all $2\leq i < j \leq k$.
\end{itemize}
In other words, as above, under an assignment $\varphi: V^{(2)} \longrightarrow \fF_q^2$, an edge $e = \{v_1 < \dots < v_k\}$ is colored red by our half-graph rule if and only if $\langle \varphi(v_iv_j) \rangle_{1 \leq i < j \leq k}$ is a half-graph pattern.

\begin{lemma}
\label{lem: half-graph-large-box}
    For every $0 < \eta < 1/6$ and all sufficiently large prime powers $q$, we have the following: if $X_1, \dots, X_{\binom{k}{2}} \subseteq \fF_q^2$ have $|X_\l| \geq M^{2 - \eta}$ for each {$\l \in [\binom{k}{2}]$}, then {$\prod_{\l=1}^{\binom{k}{2}} X_\l$} contains a half-graph pattern.
\end{lemma}

\begin{proof}
    We again apply Lemma~\ref{lem: large-box}, this time with $r=\binom{k}{2}$ and $d=2$, and $G$ now defined on $[\binom{k}{2}] \times [2]$ as follows: the rows are indexed by pairs $i < j \in [k]$, and we identify the second coordinate set $[2]$ with the pair by $1 \mapsto i$ and {$2 \mapsto j$}; then for each $i \in [k]$, we let $a_i = \min([k] \setminus \{i\})$ and take the directed out-star
    \[
        (\{a_i, i\}, i) \longrightarrow (\{i,j\}, i) \qquad \text{ for } j \in [k] \setminus \{a_i,i\}.
    \]
\end{proof}

\begin{proof}[Proofs of Theorems~\ref{thm: 3-edges-k-unif} and~\ref{thm: half-graph}]
    The probabilistic argument is the same in both cases. Put
    $m=n/k$ and let $M=|\fF_q|$.
    Since $km\leq n$, it suffices to construct a coloring with no blue
    $K_{km}^{(k)}$.
    We apply Lemma~\ref{lem: assignment-counting} with the
    following choices:
    \begin{itemize}
        \item for Theorem~\ref{thm: 3-edges-k-unif}, take
        $d=k-1$, $\scF=[k]^{(k-1)}$, and let $\scP$ be
        the set of red patterns defined by~\eqref{eq: red-pattern};
        \item for Theorem~\ref{thm: half-graph}, take
        $d=2$, $\scF=[k]^{(2)}$, and let $\scP$ be
        the set of half-graph patterns.
    \end{itemize}
    In either case, put $h=|\scF|$, $\Omega=\fF_q^d$,
    $T=M^d$, and
    \[
        \eta=\frac{1}{7d},\qquad
        b=\frac{1}{2h+1},\qquad a=\frac{1}{2h}.
    \]
    Choose a prime power $q$ so that
    \[
        m^{1/(2h+1)}\leq M^d\leq m^{1/(2h)};
    \]
    this is possible for sufficiently large $m$, for example by
    Bertrand's postulate. These choices give
    $m^b\leq T\leq m^a$ and $a<1/h$, as required.
    Moreover,
    \[
        T^{1-\eta}=M^{d-1/7},
    \]
    so the hypothesis~$(\star)$ of Lemma~\ref{lem: assignment-counting}
    follows from Lemma~\ref{lem: 3-edge-large-box} in the first
    case and Lemma~\ref{lem: half-graph-large-box} in the second.

    Now let
    \[
        N=T^{\eta m^{d-1}/(8k)}
    \]
    and choose an assignment $\varphi:[N]^{(d)}\to\Omega$
    uniformly at random. Color each $k$-set according to the
    corresponding red pattern rule above. By
    Lemma~\ref{lem: forbidden-red-config}, every $(k+1)$-set
    spans at most two red edges in the first case, while
    Lemma~\ref{lem: forbidden-halfgraph} guarantees that there is
    no red copy of $B^{(k)}$ in the second.

    On the other hand, fix a set $V\subseteq[N]$ of $km$ vertices,
    and partition it into $k$ consecutive blocks
    $V_1<\cdots<V_k$, each of size $m$. Write
    \[
        \scX=\bigcup_{F\in\scF}\prod_{i\in F}V_i.
    \]
    If $V$ is an all-blue clique, then the restriction
    $\varphi|_{\scX}$ is a $\scP$-free assignment. By
    Lemma~\ref{lem: assignment-counting}, there are at most
    $T^{hm^d-\eta m^d/4}$ such assignments, out of
    $T^{hm^d}$ assignments in total. Thus
    \[
        \bbP[V\text{ is an all-blue clique}]
        \leq T^{-\eta m^d/4}.
    \]
    Since there are at most $N^{km}$ choices for $V$,
    a union bound gives that the probability of an all-blue
    clique of size $km$ is at most
    \[
        N^{km}T^{-\eta m^d/4}
        \leq T^{\eta m^d/8-\eta m^d/4}<1.
    \]
    Consequently, there is a coloring with no forbidden red
    configuration and no blue $K_{km}^{(k)}$, giving
    \[
        r_k\bigl(k+1,3;km\bigr)>N
        \quad\text{or}\quad
        r_k\bigl(B^{(k)},K_{km}^{(k)}\bigr)>N,
    \]
    respectively. Finally, since $T=m^{\Theta_k(1)}$
    and $m=\Theta(n)$, the first choice $d=k-1$ yields
    \[
        r_k(k+1,3;n)\geq2^{\Omega_k(n^{k-2}\log n)},
    \]
    while the second choice $d=2$ yields
    \[
        r_k(B^{(k)},K_n^{(k)})\geq2^{\Omega_k(n\log n)}.
    \]
\end{proof}


\section{\texorpdfstring{$r_k(k+1, t; n)$ for $t \geq 4$}{rk(k+1, t; n) for t >= 4} and the proof of Theorem~\ref{thm: t-edges-k-unif}}
\label{sec: proofs-larger-t}

In this section, we use a stepping-up lemma of Mubayi and Suk~\cite{MS20a} along with Theorem~\ref{thm: 3-edges-k-unif} to prove Theorem~\ref{thm: t-edges-k-unif}.
This lemma reads as follows.

\begin{theorem}[Mubayi--Suk~\cite{MS20a}, Theorem 7]
\label{thm: stepping-up}
    Let $k \geq 6$ and $t \geq 4$.
    Unless $t=4$ and $k$ is odd, we have the inequality
    \begin{equation}
    \label{eq: stepping-up}
        r_k(k+1, t; 2kn) > 2^{r_{k-1}(k, t-1; n)-1}.
    \end{equation}
\end{theorem}

Now a simple iteration of this result, starting from Theorem~\ref{thm: 3-edges-k-unif}, is enough to prove Theorem~\ref{thm: t-edges-k-unif}.
We give the details below.

\begin{proof}[Proof of Theorem~\ref{thm: t-edges-k-unif}]
    Fix $k \geq 6$ and $4 \leq t \leq k-2$, and first suppose that $k-t$ is even.
    Set $j := k-t+3$.
    Starting with the bound
    \[
        r_j(j+1, 3; C_1 n) \geq 2^{C_1' n^{j-2} \log n}
        = 2^{C_1' n^{k-t+1} \log n},
    \]
    from Theorem~\ref{thm: 3-edges-k-unif}, where $C_1, C_1' > 0$ are appropriately chosen constants depending only on $k$ and $t$, we can then iterate Theorem~\ref{thm: stepping-up} $t-3$ times to obtain consecutive bounds going from
    \[
        r_j(j+1, 3; C_1 n) \to r_{j+1}(j+2, 4; C_2 n) \to \dots \to r_k(k+1, t; C_{t-2} n),
    \]
    picking up an exponential each time.
    By choosing the constants $C_1, C_1', \dots, C_{t-2}, C_{t-2}'$ to depend on $t$ and $k$ appropriately, we can ensure $C_{t-2} = 1$ so that the final bound is
    \[
        r_k(k+1, t; n) \geq T_{t-1}( C_{t-2}' n^{k-t+1} \log n)
    \]
    as promised.
    The only exceptional case is the first step, i.e., going from $r_j(j+1, 3; C_1 n) \to r_{j+1}(j+2, 4; C_2 n)$, as here the number of edges in the forbidden configuration is $4$; however, note that the uniformity $j+1 = k-t+4$ is at least $6$ by our assumption on $t$ and is even since we are in the case where $k-t$ is even, and so Theorem~\ref{thm: stepping-up} applies without issue.

    Now suppose on the other hand that $k-t$ is odd.
    If we wanted to iterate as above, here we would run into trouble at the exceptional case with $4$ edges in the forbidden configuration, as now the uniformity $k-t+4$ would be odd and Theorem~\ref{thm: stepping-up} would not apply.
    To remedy this, we use~\cite[Theorem 8]{MS20a}, which states that for odd $k > 6$,
    \begin{equation}
    \label{eq: mubayi-suk-exception}
        r_k(k+1, 4; 4n^2) > 2^{r_{k-1}(k, 3; n)-1}.
    \end{equation}
        
    Therefore, we again set $j := k-t+3$ and start with the bound
    \[
        r_j(j+1, 3; C_1 n) \geq 2^{C_1' n^{j-2} \log n}
        = 2^{C_1' n^{k-t+1} \log n},
    \]
    from Theorem~\ref{thm: 3-edges-k-unif}, but then use~\eqref{eq: mubayi-suk-exception} to get
    \[
        r_{j+1}(j+2, 4; C_2 n) > T_3( C_2' n^{\frac{k-t+1}{2}} \log n);
    \]
    indeed, now the uniformity $j+1 = k-t+4$ is greater than $6$ (since $k-t \geq 2$ is odd), and the constant $C_2'$ accounts for the change $\log n = \Theta_{k, t}(\log (\sqrt{C_2 n}/2))$.
    With this taken care of, we can then iterate Theorem~\ref{thm: stepping-up} $t-4$ more times as above to obtain successive bounds going from 
    \[
        r_{j+1}(j+2, 4; C_2 n) \to r_{j+2}(j+3, 5; C_3n) \to \dots \to r_k(k+1, t; C_{t-2} n),
    \]
    picking up an exponential each time.
    By choosing the constants depending on $k$ and $t$ appropriately once again, we can take $C_{t-2} = 1$ so that the final bound is
    \[
        r_k(k+1, t; n) \geq T_{t-1}( C_{t-2}' n^{\frac{k-t+1}{2}} \log n).
    \]
    This completes the proof.
\end{proof}


\section{Acknowledgements}
\label{sec: conclusion}

We would like to thank David Conlon for helpful discussions, and Anurag Bishnoi for directing our attention to the references~\cite{BMM20, Ott, Smith} for more background on the directed graph $D_q$.

We also acknowledge the use of ChatGPT (OpenAI) in assisting with the proofs of Theorems~\ref{thm: 3-edges-k-unif} and~\ref{thm: half-graph}. 
The main idea of combining the digraph behind Brada\v{c}'s recent construction from \cite{Brad26+} with the classical Erd\H{o}s--Hajnal ideas from~\cite{EH72} in order to obtain a simpler proof of Theorem~\ref{thm: k4-} is due to the authors.

The initial draft was written from scratch by the authors once they had digested the proofs of Theorems~\ref{thm: 3-edges-k-unif} and~\ref{thm: half-graph}, after which ChatGPT was used to check for typos and suggest simplifications. 
The authors take full responsibility for the contents of the paper.

\subsection*{Note added}
When finalizing the initial draft of this manuscript, we learned that similar results were obtained independently by Bai, Hu, Liu, and Wang~\cite{BHLW++}.
In particular, they have also proved Theorems~\ref{thm: 3-edges-k-unif},~\ref{thm: half-graph}, and~\ref{thm: t-edges-k-unif}, along with a new stepping-up lemma which resolves the case when $k-t$ is odd.


\end{document}